\documentclass[12pt]{article}
\usepackage{amsmath,amssymb}
\newtheorem{theorem}{\hspace*{\parindent}Theorem}
\newtheorem{lemma}{\hspace*{\parindent}Lemma}
\newtheorem{corollary}{\hspace*{\parindent}Corollary}
\newtheorem{conjecture}{\hspace*{\parindent}Conjecture}
\def\N{\mathbb{N}}
\title{Normalized incomplete beta function: log-concavity in parameters and other properties}
\author{D.B.\:Karp$^{\rm a,b}$\footnote{E-mail: D. Karp -- \emph{dimkrp@gmail.com}}
\\[10pt]\small{\textit{$\phantom{1}^a$Far Eastern Federal University, 8 Sukhanova street, Vladivostok, 690950, Russia}}\\\small{\textit{$\phantom{1}^b$Institute of Applied Mathematics, FEBRAS, 7 Radio Street, Vladivostok,  690041, Russia}}}
\date{}
\begin{document}
\maketitle

\bigskip

\begin{center}
\parbox{12cm}{
\small\textbf{Abstract.} The normalized incomplete beta function can be defined either as cumulative distribution
function of beta density or as the Gauss hypergeometric function with one of the upper parameters equal to unity.
Logarithmic concavity/convexity of this function in parameters was established by Finner and Roters in 1997. Their proof
is indirect and rather difficult; it is based on generalized reproductive property of certain more general distributions.
These authors remark that these results ''seems to be very hard to obtain by usual analytic methods''.  In the first part
of this paper we provide such proof based on standard tools of analysis.  In the second part we go one step further
and investigate the sign of generalized Tur\'{a}n determinants formed by shifts of the normalized incomplete beta function.
Under some additional restrictions we demonstrate that these coefficients are of the same sign.  We further conjecture that
such restrictions can be removed without altering the results.  Our method of proof also leads to various companion results
which may be of independent interest.  In particular, we establish linearization formulas and two-sided bounds for the above
mentioned Tur\'{a}n determinants.  Further, we find two combinatorial style identities for finite sums which we believe
to be new.}
\end{center}

\bigskip


\bigskip


\bigskip

\paragraph{1. Motivation and introduction.}
The beta distribution is perhaps the single most important continuous compactly supported probability distributions.
Its particular cases include the Wigner semicircle, the Marchenko-Pastur and the arcsine laws; it is important in
Bayesian analysis as conjugate prior to binomial and geometric distributions \cite{BayesBook};
it plays a role in a large number of applications ranging from population genetics to project management.
The beta distribution is defined by the density $x^{a-1}(1-x)^{b-1}/B(a,b)$ supported on $[0,1]$, where $B(a,b)$
is Euler's beta function \cite{AAR}. The cumulative distribution function (CDF) of the beta density is given by
the normalized incomplete beta function
$$
I_x(a,b)=\frac{B_{x}(a,b)}{B(a,b)}=\frac{\int_0^xt^{a-1}(1-t)^{b-1}dt}{\int_0^1t^{a-1}(1-t)^{b-1}dt}.
$$
In a recent paper \cite{Garibaldi} containing the mathematical analysis supporting a recent investigation related to
lottery frauds in Florida, the authors solve the following optimization problem
$$
\sum_{i}^{m}c_i\alpha_i\to\min~~\text{under constraint}~
f(\alpha_1,\alpha_2,\ldots,\alpha_m)=\prod\limits_{i=1}^{m}I_{p_i}(w_i,\alpha_i-w_i+1)\ge\varepsilon.
$$
For efficient numerical solution of such problem it is desirable that the function $f$ be quasi-concave
which implies that the feasible set of the problem is convex and local minimizer is also global.
Since monotonic transformations do not alter quasi-concavity it suffices then to show that $\log{f}$ is a
concave function.  This fact, in turn, is implied by log-concavity of $b\to{I_{x}(a,b)}$ on $(0,\infty)$.
The purpose of this paper is to investigate the log-concavity/convexity properties of the function $I_{x}(a,b)$
viewed as the function of parameters $a$ and $b$. In particular, we demonstrate that $b\to{I_{x}(a,b)}$ is
indeed log-concave on $(0,\infty)$, while $a\to{I_{x}(a,b)}$ is log-convex on $(0,\infty)$ for $b\in(0,1)$ and
log-concave for $b>1$. Log-concavity of  $b\to{I_{x}(a,b)}$ is equivalent to the positivity of
the generalized Tur\'{a}n determinant
$$
I_x(a,b+\alpha)I_x(a,b+\beta)-I_x(a,b)I_x(a,b+\alpha+\beta)
$$
for all $\alpha,\beta>0$, while log-convexity (log-concavity) of $a\to{I_{x}(a,b)}$ is equivalent to negativity (positivity) of
$$
I_x(a+\alpha,b)I_x(a+\beta,b)-I_x(a,b)I_x(a+\alpha+\beta,b).
$$
We go one step further and study the signs of the power series coefficients (in powers of $x$) of the above Tur\'{a}n determinants.
Under some additional restrictions we demonstrate that these coefficients are of the same sign.  We further conjecture
that such restrictions can be removed without altering the results. Our method of proof also leads to various companion results
which may be of independent interest.  In particular, we establish linearization formulas and two-sided bounds for the above
Tur\'{a}n determinants.  Further, we find two combinatorial style identities for finite sums which we believe
to be new.

When this paper was nearly completed we discovered that log-concavity/convexity of the CDF of beta distribution has been demonstrated
in 1997 by Finner and Roters in their fundamental work \cite{FR97}. Partial results in this direction have been previously given
by the same authors in \cite{FR93} and by Das Gupta and Sarkar in \cite{Gupta-Sarkar}.  However, their log-concavity
proofs are very involved and indirect - they appear as a by-product of generalized reproductive property of certain
probability measures. Finner and Roters  note in their paper that ''we are now able to conclude a result for
the Beta distribution which seems to be very hard to obtain by usual analytic methods''.  The proofs presented in
this paper are precisely ''by usual analytic methods''. The power series  coefficients of the Tur\'{a}n determinants
have not been considered by the above authors, so our results in this direction strengthen the achievements of \cite{FR97}.

We conclude the introduction by presenting an alternative  expression for the normalized incomplete beta function
in terms of the  Gauss hypergeometric function \cite{AAR}
$$
{_2F_1}(\alpha,\beta;\gamma;z)=\sum\limits_{n=0}^{\infty}\frac{(\alpha)_n(\beta)_n}{(\gamma)_nn!}z^n.
$$
Indeed, using Euler-Pochhammer integral representation
$$
{_2F_1}(\alpha,\beta;\gamma;x)=\frac{\Gamma(\gamma)}{\Gamma(\beta)\Gamma(\gamma-\beta)}\int_0^1u^{\beta-1}(1-u)^{\gamma-\beta-1}(1-ux)^{-\alpha}du
$$
we get by a change of variable and application of Euler's formula $B(a,b)=\Gamma(a)\Gamma(b)/\Gamma(a+b)$:
$$
I_x(a,b)=\frac{\Gamma(a+b)x^{a}}{\Gamma(a+1)\Gamma(b)}{_2F_1}(1-b,a;a+1;x).
$$
Further, applying another Euler's  formula
$$
{_2F_1}(\alpha,\beta;\gamma;x)=(1-x)^{\gamma-\alpha-\beta}{_2F_1}(\gamma-\alpha,\gamma-\beta;\gamma;x)
$$
we obtain the well-known representation
\begin{equation}\label{eq:I2F1}
I_x(a,b)=\frac{\Gamma(a+b)x^{a}}{\Gamma(a+1)\Gamma(b)}(1-x)^{b}{_2F_1}(a+b,1;a+1;x)
=x^{a}(1-x)^{b}\sum\limits_{n=0}^{\infty}\frac{\Gamma(a+b+n)}{\Gamma(b)\Gamma(a+1+n)}x^n.
\end{equation}

\paragraph{2. Log-concavity of $I_x(a,b)$ in $a$ and $b$.}
The proof of the next theorem has been inspired by the log-concavity proof for incomplete gamma function given in \cite{Alzer98}. See a related result in \cite{AB2012}.
\begin{theorem}\label{th:logconc-b}
For each fixed $a>0$ and $x\in(0,1)$ the function $b\to{I_x(a,b)}$ is strictly log-concave on $(0,\infty)$.
\end{theorem}
\textbf{Proof.}  We need to show that
$$
\frac{\partial^2}{\partial{b}^2}\log(I_x(a,b))<0.
$$
It is easy to compute
\begin{equation}\label{eq:logderivB}
\frac{\partial^2}{\partial{b}^2}\log(B(a,b))=\psi'(b)-\psi'(a+b)>0,
\end{equation}
where $\psi(z)=\Gamma'(z)/\Gamma(z)$ is logarithmic derivative of the gamma function and the inequality follows from
the fact that $z\to\psi'(z)$ is decreasing on $(0,\infty)$ according to the representation
$$
\psi'(z)=\sum\limits_{n=0}^{\infty}\frac{1}{(n+z)^2}.
$$
Further,
\begin{multline*}
B_{x}(a,b)^2\frac{\partial^2}{\partial{b}^2}\log(I_x(a,b))=B_{x}(a,b)\int_0^xt^{a-1}(1-t)^{b-1}[\log(1-t)]^2dt-
\\
\left[\int_0^xt^{a-1}(1-t)^{b-1}\log(1-t)dt\right]^2-B_{x}(a,b)^2(\psi'(b)-\psi'(a+b))=:U_{a,b}(x)
\end{multline*}
We have $U_{a,b}(0)=U_{a,b}(1)=0$.  The second equality can be seen by setting $x=1$ and then differentiating which is legitimate by continuous double differentiability or by setting $x=1$ in the above formula and computing the integrals.
We aim to demonstrate that $U_{a,b}(x)$ is decreasing on $(0,x_1)$ and increasing on $(x_1,1)$ for some $x_1\in(0,1)$.  This will prove that $U_{a,b}(x)<0$ for $x\in(0,1)$ for each $a,b>0$. We have
\begin{multline*}
V_{a,b}(x):=\frac{U_{a,b}'(x)}{x^{a-1}(1-x)^{b-1}}=\int_0^xt^{a-1}(1-t)^{b-1}[\log(1-t)]^2dt+B_{x}(a,b)[\log(1-x)]^2
\\
-2\log(1-x)\int_0^xt^{a-1}(1-t)^{b-1}\log(1-t)dt-2B_{x}(a,b)(\psi'(b)-\psi'(a+b)).
\end{multline*}
It is rather straightforward to see that $V_{a,b}(0)=0$ and $V_{a,b}(1)=+\infty$.   We aim to demonstrate that $V_{a,b}(x)$
 is decreasing on $(0,x_2)$ and increasing on $(x_2,1)$ for some $x_2\in(0,1)$. This will imply that it
 changes sign on $(0,1)$ exactly once from minus to plus which in turns implies the same conclusion for $U_{a,b}'(x)$.
 Taking the next derivative we obtain after cancelations:
\begin{multline*}
W_{a,b}(x):=\frac{1}{2}(1-x)V_{a,b}'(x)=-B_x(a,b)\log(1-x)
\\
+\int_0^xt^{a-1}(1-t)^{b-1}\log(1-t)dt-x^{a-1}(1-x)^b(\psi'(b)-\psi'(a+b)).
\end{multline*}
It is clear that  $W_{a,b}(1)=+\infty$.  The value at zero depends on $a$: if $0<a<1$ then $W_{a,b}(0)=-\infty$; if $a=1$ then $W_{a,b}(0)=-(\psi'(b)-\psi'(a+b))<0$; if $a>1$ then $W_{a,b}(0)=0$. We will demonstrate below that $W_{a,b}(x)$ is increasing on $(0,1)$ when $0<a\leq{1}$ and $W_{a,b}(x)$ is decreasing on $(0,x_3)$ and increasing on $(x_3,1)$ for some $x_3\in(0,1)$ when $a>1$. This will imply that it changes sign on $(0,1)$ exactly once from minus to plus which in turns implies the same conclusion for $V_{a,b}'(x)$. Taking the next derivative we get:
$$
Z_{a,b}(x):=(1-x)W_{a,b}'(x)=B_x(a,b)-(\psi'(b)-\psi'(a+b))x^{a-2}(1-x)^{b}((a-1)(1-x)-bx).
$$
We have $Z_{a,b}(1)=B(a,b)>0$. At $x=0$ five cases reveal themselves:

Case~I: if $0<a<1$ then $Z_{a,b}(0)=+\infty$.  In this case $Z_{a,b}(x)>0$ on $(0,1)$ since $((a-1)(1-x)-bx)<0$ and $\psi'(b)-\psi'(a+b)>0$;

Case~II: if  $a=1$ then  $Z_{a,b}(0)=b(\psi'(b)-\psi'(a+b))>0$;  We will show below that  $Z_{a,b}'(x)=0$ on $(0,1)$ so that $Z_{a,b}(x)=b(\psi'(b)-\psi'(a+b))>0$ for all $x\in(0,1)$;

Case~III:  if $1<a<2$ then $Z_{a,b}(0)=-\infty$;

Case~IV: if $a=2$ then $Z_{a,b}(0)=-(\psi'(b)-\psi'(a+b))<0$;

Case~V:  if $a>2$ then $Z_{a,b}(0)=0$;

Taking one more derivative we get:
\begin{multline*}
Q_{a,b}(x):=Z_{a,b}'(x)x^{3-a}(1-x)^{1-b}
\\
=x^2-(\psi'(b)-\psi'(a+b))((a+b-1)^2x^2-(a-1)(2a+2b-3)x+(a-1)(a-2)).
\end{multline*}
Straightforward calculation yields:
$$
Q_{a,b}(0)=-(\psi'(b)-\psi'(a+b))(a-1)(a-2),~~~~Q_{a,b}(1)=1-b^2(\psi'(b)-\psi'(a+b)).
$$
Since $a\to{Q_{a,b}(1)}$ is decreasing on $(0,\infty)$, $Q_{0,b}(1)=1$ and $Q_{1,b}(1)=0$ because $\psi'(1+b)=\psi'(b)-1/b^2$,
these formulas lead to the following conclusions.

Case~I: if $0<a<1$ then $Q_{a,b}(0)<0$ and $Q_{a,b}(1)>0$ for all $b>0$.  Since $Q_{a,b}(x)$ is a quadratic this implies that it has exactly one change of sign on $(0,1)$.

Case II: if $a=1$  then $Q_{a,b}(x)=0$ on $[0,1]$ for all $b>0$ since $\psi'(1+b)=\psi'(b)-1/b^2$.

Case III: if $1<a<2$ then $Q_{a,b}(0)>0$ and $Q_{a,b}(1)<0$ for all $b>0$.  Since $Q_{a,b}(x)$ is a quadratic this implies that it has exactly one change of sign on $(0,1)$.

Case IV: If $a=2$ then $Q_{a,b}(0)=0$ and $Q_{a,b}(1)<0$ for all $b>0$. It follows that $Q_{a,b}(x)>0$ on $(0,\alpha)$ and $Q_{a,b}(x)<0$ on $(\alpha,1)$ for some $\alpha\in(0,1)$.  An apparent alternative $Q_{a,b}(x)<0$ on $(0,1)$ cannot hold since $Z_{a,b}(x)$ would then  be decreasing on $(0,1)$ while $Z_{a,b}(0)<0$ and $Z_{a,b}(1)>0$ rules out such possibility;

Case V: $a>2$  then $Q_{a,b}(0)<0$ and $Q_{a,b}(1)<0$ for all $b>0$.  This implies that $Q_{a,b}(x)$ follows the sign pattern $(-,+,-)$ on $(0,1)$.  An apparent alternative $Q_{a,b}(x)<0$ on $(0,1)$ cannot hold since $Z_{a,b}(x)$ would  then be decreasing on $(0,1)$ while  $Z_{a,b}(0)=0$ and  $Z_{a,b}(1)>0$ rules out such possibility;

Looking carefully at each case we see that for all $b>0$ the function $Z_{a,b}(x)$ is positive on $(0,1)$ if $0<a\leq{1}$
and changes sign from minus to plus if $a>1$.  This implies that $W_{a,b}'(x)$ behaves in the same way, so
that in all cases $W_{a,b}(x)$ changes sign exactly once from minus to plus. This in turn implies that $V_{a,b}'(x)$
has the same behavior and hence  $V_{a,b}(x)$ is first decreasing and then increasing. In view of the boundary
values at $x=0$ and $x=1$ this means that $V_{a,b}(x)$ changes sign exactly once from minus to plus.  Finally, it
follows that $U_{a,b}'(x)$ is negative on $(0,x_1)$ and positive on $(x_1,1)$ for some $x_1\in(0,1)$,
yielding $U_{a,b}(x)<0$ on $(0,1)$ as claimed.~~$\hfill\square$

Definition of $I_x(a,b)$  immediately leads to the reflection formula
$$
I_x(a,b)=1-I_{1-x}(b,a).
$$
This implies
\begin{multline*}
\frac{\partial^2}{\partial{a}^2}\log(I_x(a,b))=[I_x(a,b)]^{-2}\biggl\{I_{1-x}(b,a)\frac{\partial^2}{\partial{a}^2}I_{1-x}(b,a)
-\left[\frac{\partial}{\partial{a}}I_{1-x}(b,a)\right]^2\biggr\}
\\
-[I_x(a,b)]^{-2}\frac{\partial^2}{\partial{a}^2}I_{1-x}(b,a).
\end{multline*}
The first term is negative according to Theorem~\ref{th:logconc-b}.
The second term, however, may change sign depending on the values of $b$ and $x$ as demonstrated
by numerical evidence.  Hence, we cannot draw any definitive conclusion about log-concavity of $a\to{I_x(a,b)}$
from the reflection formula.  Indeed, it turns out that the result depends on the value of $b$.
Nevertheless, a method similar to that used in the proof of Theorem~\ref{th:logconc-b} also works here combined with the
following lemma.
\begin{lemma}\label{lm:psi}
For fixed $\alpha>0$ define $f_{\alpha}(x):(0,\infty)\to(-\infty,\infty)$ by
$$
f_{\alpha}(x)=(\psi(x+\alpha)-\psi(x))^2+\psi'(x+\alpha)-\psi'(x).
$$
Then $f_{\alpha}(x)<0$ if $0<\alpha<1$ and $f_{\alpha}(x)>0$ if $\alpha>1$.
\end{lemma}
\textbf{Proof.} Using recurrence relations $\psi(x+1)=\psi(x)+1/x$ and $\psi(x+1)=\psi(x)-1/x^2$, we
calculate
\begin{multline*}
f(x+1)-f(x)=
\\
(\psi(x+\alpha+1)-\psi(x+1))^2-(\psi(x+\alpha)-\psi(x))^2+\psi'(x+\alpha+1)-\psi'(x+\alpha)+\psi'(x)-\psi'(x+1)
\\
=(\psi(x+\alpha+1)-\psi(x+1)+\psi(x+\alpha)-\psi(x))(\psi(x+\alpha+1)-\psi(x+1)-\psi(x+\alpha)+\psi(x))
\\
+\frac{1}{x^2}-\frac{1}{(x+\alpha)^2}
=\left(2\psi(x+\alpha)+\frac{1}{x+\alpha}-2\psi(x)-\frac{1}{x}\right)\left(\frac{1}{x+\alpha}-\frac{1}{x}\right)+\frac{1}{x^2}-\frac{1}{(x+\alpha)^2}
\\
=\frac{\alpha}{x(x+\alpha)}\left(2\psi(x)-2\psi(x+\alpha)+\frac{\alpha}{x(x+\alpha)}\right)+\frac{\alpha(2x+\alpha)}{x^2(x+\alpha)^2}
\\
=\frac{2\alpha}{x(x+\alpha)}
\left(\psi(x)-\psi(x+\alpha)+\frac{1}{x}\right).
\end{multline*}
Writing $g_{\alpha}(x)$ for the function in parentheses and utilizing the integral representations
$$
\psi(x)=-\gamma+\int_0^{\infty}\frac{e^{-t}-e^{-xt}}{1-e^{-t}}dt~~~\text{and}~~~~\frac{1}{x}=\int_0^{\infty}e^{-xt}dt,
$$
where $\gamma$ stands for Euler-Mascheroni constant, we get:
$$
g_{\alpha}(x)=\int_0^{\infty}\frac{e^{-xt}(e^{-\alpha{t}}-e^{-t})}{1-e^{-t}}dt.
$$
The last formula makes it obvious that $g_{\alpha}(x)>0$ for $0<\alpha<1$ and $g_{\alpha}(x)<0$ for $\alpha>1$.
Further, from the asymptotic formulas
$$
\psi(x)=\log(x)-\frac{1}{2x}+O(x^{-2}),~~~~~\psi'(x)=\frac{1}{x}+\frac{1}{2x^2}+O(x^{-3}),~~~x\to\infty,
$$
we conclude that $\lim\limits_{x\to\infty}f(x)=0$.  Altogether this implies for $0<\alpha<1$:
$$
f(x)<f(x+1)<f(x+2)<\cdots<\lim\limits_{m\to\infty}f(x+m)=0,
$$
and for $\alpha>1$:
$$
f(x)>f(x+1)>f(x+2)>\cdots>\lim\limits_{m\to\infty}f(x+m)=0,
$$
which proves the lemma.~~$\hfill\square$

Stronger results for a function similar to $f_{\alpha}(x)$ but still slightly different from it
can be found in \cite[Theorem~1.2]{Qi-Guo}.

\begin{theorem}\label{th:logconc-a}
Suppose $x\in(0,1)$ is fixed. Then:

\emph{(a)} if $0<b<1$ the function $a\to{I_x(a,b)}$ is strictly log-convex on $(0,\infty)$.

\emph{(b)} if $b=1$, the function $I_x(a,b)=x^a$ and so is log-neutral.

\emph{(c)} if $b>1$  the function $a\to{I_x(a,b)}$ is strictly log-concave on $(0,\infty)$.
\end{theorem}
\textbf{Proof.}  The case $b=1$ can be verified directly.  It remains to show that
$$
\frac{\partial^2}{\partial{a}^2}\log(I_x(a,b))>0~~\text{if}~0<b<1~~\text{and}~~\frac{\partial^2}{\partial{a}^2}\log(I_x(a,b))<0
~~\text{if}~b>1.
$$
For completeness, however, we will include $b=1$ into the forgoing considerations.  By symmetry we have from (\ref{eq:logderivB})
$$
\frac{\partial^2}{\partial{a}^2}\log(B(a,b))=\psi'(a)-\psi'(a+b)>0.
$$
Further,
\begin{multline*}
B_{x}(a,b)^2\frac{\partial^2}{\partial{a}^2}\log(I_x(a,b))=B_{x}(a,b)\int_0^xt^{a-1}(1-t)^{b-1}[\log(t)]^2dt-
\\
\left[\int_0^xt^{a-1}(1-t)^{b-1}\log(t)dt\right]^2-B_{x}(a,b)^2(\psi'(a)-\psi'(a+b))=:U_{a,b}(x).
\end{multline*}
 Direct verification yields $U_{a,b}(0)=U_{a,b}(1)=0$.  Next, we have
\begin{multline*}
V_{a,b}(x):=\frac{U_{a,b}'(x)}{x^{a-1}(1-x)^{b-1}}=\int_0^xt^{a-1}(1-t)^{b-1}[\log(t)]^2dt+B_{x}(a,b)[\log(x)]^2
\\
-2\log(x)\int_0^xt^{a-1}(1-t)^{b-1}\log(t)dt-2B_{x}(a,b)(\psi'(a)-\psi'(a+b)).
\end{multline*}
Repeatedly  using L'Hopital's rule we compute the limits:
\begin{multline*}
\lim\limits_{x\to0}\frac{B_{x}(a,b)}{[\log(x)]^{-2}}=\lim\limits_{x\to0}\frac{x^{a}(1-x)^{b-1}}{2(\log(1/x))^{-3}}
=\frac{1}{2}\lim\limits_{x\to0}\frac{(\log(1/x))^{3}}{x^{-a}}
\\
=\frac{1}{2}\lim\limits_{x\to0}\frac{3(\log(1/x))^{2}}{ax^{-a}}
=\frac{3}{2}\lim\limits_{x\to0}\frac{2\log(1/x)}{a^2x^{-a}}=3\lim\limits_{x\to0}\frac{1}{a^3x^{-a}}=0.
\end{multline*}
Similarly,
$$
\lim\limits_{x\to0}\log(x)\int_0^xt^{a-1}(1-t)^{b-1}\log(t)dt=0,
$$
so that $V_{a,b}(0)=0$.  Further,
$$
\int_0^1t^{a-1}(1-t)^{b-1}[\log(t)]^2dt=B(a,b)[(\psi(a)-\psi(a+b))^2+\psi'(a)-\psi'(a+b)],
$$
so that
$$
V_{a,b}(1)=B(a,b)[(\psi(a)-\psi(a+b))^2-(\psi'(a)-\psi'(a+b))].
$$
It follows from Lemma~\ref{lm:psi} that (a) $V_{a,b}(1)<0$ for $0<b<1$; (b) $V_{a,b}(1)=0$ for $b=1$; (c) $V_{a,b}(1)>0$ for $b>1$;
 We aim to demonstrate that $V_{a,b}(x)$ has precisely one change of sign in cases (a) and (c) and is identically zero in case (b).
 Taking the next derivative we obtain after cancelations:
\begin{multline*}
W_{a,b}(x):=\frac{x}{2}V_{a,b}'(x)=B_x(a,b)\log(x)
\\
-\int_0^xt^{a-1}(1-t)^{b-1}\log(t)dt-x^{a}(1-x)^{b-1}(\psi'(a)-\psi'(a+b)).
\end{multline*}
The boundary values are:  $W_{a,b}(0)=0$ and
$$
W_{a,b}(1)=\left\{\!\!\begin{array}{ll}
-\infty, & 0<b<1
\\
0, & b=1,
\\
B(a,b)(\psi(a+b)-\psi(a))>0, & b>1.
\end{array}
\right.
$$
These values follow from the evaluation
$$
\int_0^xt^{a-1}(1-t)^{b-1}\log(t)dt=B(a,b)(\psi(a)-\psi(a+b))
$$
and the recurrence relations $\psi(a+1)-\psi(a)=1/a$, $\psi'(a+1)-\psi'(a)=-1/a^2$.
Taking the next derivative we get
$$
Z_{a,b}(x)=xW'_{a,b}(x):=B_x(a,b)+(\psi'(a)-\psi'(a+b))x^{a}(1-x)^{b-2}((b-1)x-a(1-x)).
$$
Clearly, $Z_{a,b}(0)=0$. At $x=1$ five cases reveal themselves:

Case~I: if $0<b<1$ then $Z_{a,b}(1)=-\infty$;

Case~II: if  $b=1$ then  $Z_{a,b}(1)=0$;  We will show below that  $Z_{a,b}'(x)=0$ on $(0,1)$ so that $Z_{a,b}(x)=0$ for all $x\in(0,1)$;

Case~III:  if $1<b<2$ then $Z_{a,b}(1)=+\infty$;

Case~IV: if $b=2$ then $Z_{a,b}(1)=B(a,b)+\psi'(a)-\psi'(a+b)>0$;

Case~V:  if $b>2$ then $Z_{a,b}(1)=B(a,b)>0$.

Finally, computing  one more derivative and making some rearrangements we arrive at
\begin{multline*}
Q_{a,b}(x)=\frac{Z_{a,b}'(x)}{x^{a-1}(1-x)^{b-3}}
\\
=(1-x)^{2}-(\psi'(a)-\psi'(a+b))((a+b-1)^2x^2-((2a+1)(b-1)+2a^2)x+a^2).
\end{multline*}
The boundary values are
$$
Q_{a,b}(0)=1-a^2(\psi'(a)-\psi'(a+b)),~~~~Q_{a,b}(1)=-(\psi'(a)-\psi'(a+b))(b-1)(b-2).
$$
Since $b\to{Q_{a,b}(0)}$ is decreasing on $(0,\infty)$, $Q_{a,0}(0)=1$ and $Q_{a,1}(0)=0$ because $\psi'(1+a)=\psi'(a)-1/a^2$,
these formulas lead to the following conclusions.

Case~I: if $0<b<1$ then $Q_{a,b}(0)>0$ and $Q_{a,b}(1)<0$ for all $a>0$.  Since $Q_{a,b}(x)$ is a quadratic this implies that it has exactly one change of sign on $(0,1)$.

Case II: if $b=1$  then $Q_{a,b}(x)=0$ on $[0,1]$ for all $a>0$ since $\psi'(1+a)=\psi'(a)-1/a^2$.

Case III: if $1<b<2$ then $Q_{a,b}(0)<0$ and $Q_{a,b}(1)>0$ for all $b>0$.  Since $Q_{a,b}(x)$ is a quadratic this implies that it has exactly one change of sign on $(0,1)$.

Case IV: If $b=2$ then $Q_{a,b}(0)<0$ and $Q_{a,b}(1)=0$ for all $a>0$. It follows that $Q_{a,b}(x)<0$ on $(0,\alpha)$
and $Q_{a,b}(x)>0$ on $(\alpha,1)$ for some $\alpha\in(0,1)$.  An apparent alternative $Q_{a,b}(x)<0$ on $(0,1)$ cannot
hold since $Z_{a,b}(x)$ would then be decreasing on $(0,1)$ while $Z_{a,b}(0)=0$ and $Z_{a,b}(1)>0$ rules out such possibility.

Case V: $b>2$  then $Q_{a,b}(0)<0$ and $Q_{a,b}(1)<0$ for all $a>0$.  This implies that $Q_{a,b}(x)$ follows the sign
pattern $(-,+,-)$ on $(0,1)$.  An apparent alternative $Q_{a,b}(x)<0$ on $(0,1)$ cannot hold since $Z_{a,b}(x)$
would then be decreasing on $(0,1)$ while $Z_{a,b}(0)=0$ and $Z_{a,b}(1)>0$ rules out such possibility.

Tracing back the primitives $Q_{a,b}\to{Z_{a,b}}\to{W_{a,b}}\to{V_{a,b}}\to{U_{a,b}}$ we  see that for all $a>0$
 the function $Z_{a,b}(x)$ changes sign from plus to minus on $(0,1)$ if $0<b<1$ is identically zero if $b=1$ and
 changes sign from minus to plus if $b>1$.  This implies that $W_{a,b}(x)$ and $V_{a,b}$(x) behave in the
 same way, so that $U_{a,b}(x)$ is first increasing  and then decreasing for $0<b<1$, is identically zero if $b=1$
 and is first decreasing and then increasing if $b>1$.  In view of the boundary values at $x=0$ and $x=1$ this
 means that $U_{a,b}(x)>0$ on $(0,1)$ for $0<b<1$ and  $U_{a,b}(x)<0$ for $b>1$ as claimed.~$\hfill\square$

\paragraph{3. Coefficient-wise log-concavity of $I_x(a,b)$.}
As mentioned in the introduction, logarithmic concavity of $b\to{I_x(a,b)}$ demonstrated
in Theorem~\ref{th:logconc-b} is equivalent to the inequality
\begin{equation}\label{eq:Wright}
\phi_{a,b,\alpha,\beta}(x):=I_x(a,b+\alpha)I_x(a,b+\beta)-I_x(a,b)I_x(a,b+\alpha+\beta)>0
\end{equation}
valid for all positive values of $a$, $b$, $\alpha$ and $\beta$ and each $x\in(0,1)$.  This inequality
expresses the strict Wright log-concavity of $b\to{I_x(a,b)}$ - a property equivalent to log-concavity
for all continuous functions \cite[Chapter~I.4]{MPF}, \cite[Section~1.1]{PPT}. On the other hand, formula (\ref{eq:I2F1})
implies that the function of $\phi_{a,b,\alpha,\beta}(x)$ defined in (\ref{eq:Wright}) has the power
series expansion
\begin{equation}\label{eq:phi-expanded}
\phi_{a,b,\alpha,\beta}(x)=x^{2a}(1-x)^{2b+\alpha+\beta}\sum\limits_{k=0}^{\infty}\phi_kx^k.
\end{equation}
It is then natural to ask whether and when the coefficients $\phi_k$ are positive.
In general, for a given formal power series
\begin{equation}\label{eq:fgeneral}
f(\mu;x)=\sum\limits_{k=0}^{\infty}f_k(\mu)x^k
\end{equation}
with nonnegative coefficients $f_k(\mu)$ which depend continuously on a parameter
$\mu$ from a real interval $I$, we will say that $\mu\to{f(\mu;x)}$ is coefficient-wise Wright log-concave
on $I$ for $\alpha\in{A}\subseteq[0,\infty)$, $\beta\in{B}\subseteq[0,\infty)$
if the formal power series for the product difference
$$
f(\mu+\alpha;x)f(\mu+\beta;x)-f(\mu;x)f(\mu+\alpha+\beta;x)
$$
has nonnegative coefficients at all powers of $x$ for $\alpha\in{A}$, $\beta\in{B}$ and all $\mu\in{I}$.
If coefficients are strictly positive we say that the corresponding property holds strictly.
If coefficients are non-positive we talk about coefficient-wise Wright log-convexity.
These concepts using slightly different terminology were introduced by the author (jointly with S.I.\:Kalmykov) in
\cite{KKJMAA}.

Conclusions of Theorem~\ref{th:logconc-a} can be restated in terms of the function
\begin{equation}\label{eq:psi-def}
\psi_{a,b,\alpha,\beta}(x):=I_x(a+\alpha,b)I_x(a+\beta,b)-I_x(a,b)I_x(a+\alpha+\beta,b)
=x^{2a+\alpha+\beta}(1-x)^{2b}\sum\limits_{k=0}^{\infty}\psi_kx^k
\end{equation}
as follows: if $a>0$ and $0<x<1$ then $\psi_{a,b,\alpha,\beta}(x)<0$
for $0<b<1$ and $\psi_{a,b,\alpha,\beta}(x)>0$ for $b>1$. Just like with $\phi_{a,b,\alpha,\beta}(x)$
it is then natural to consider the coefficient-wise log-convexity and log-concavity of $a\to x^{-a}(1-x)^{-b}I_x(a,b)$,
i.e. study the sign of the coefficients $\psi_k$.

The main purpose of this section is to give partial solution to the following conjectures
representing substantial strengthenings of Theorems~\ref{th:logconc-b} and \ref{th:logconc-a},
respectively.

\begin{conjecture}\label{cj:phi}
The coefficients $\phi_k$ are positive for all $a,b,\alpha,\beta>0$ and $x\in(0,1)$, so that
$b\to x^{-a}(1-x)^{-b}I_x(a,b)$ is strictly coefficient-wise Wright log-concave on $(0,\infty)$
for these values of parameters.
\end{conjecture}

\begin{conjecture}\label{cj:psi}
The coefficients $\psi_k$ are negative for all $a,\alpha,\beta>0$ and $x\in(0,1)$ if $0<b<1$ and
positive if $b>1$, so that $a\to x^{-a}(1-x)^{-b}I_x(a,b)$ is strictly coefficient-wise Wright log-convex
on $(0,\infty)$ for $\alpha,\beta>0$, $x\in(0,1)$ and $0<b<1$; and $a\to x^{-a}(1-x)^{-b}I_x(a,b)$ is
strictly coefficient-wise Wright log-concave on $(0,\infty)$ for $\alpha,\beta>0$, $x\in(0,1)$ and $b>1$.
\end{conjecture}

Our approach also leads to a number of related results which may be of independent interest.
These include alternative formulas for the coefficients $\phi_k$ and $\psi_k$, linearization identities
for the product differences of  hypergeometric functions, two-sided inequalities for the Tur\'{a}n
determinants formed by normalized incomplete beta functions and two presumably new combinatorial identities.

We will need the following lemma (see \cite[Lemma~3]{KKJMAA}).
\begin{lemma}\label{lm:q-log}
Suppose $\mu\to{f(\mu;x)}$ is  coefficient-wise Wright log-concave on $[0,\infty)$ for $\alpha=1$, $\beta\ge0$, i.e.
$$
f(\mu+1;x)f(\mu+\beta;x)-f(\mu;x)f(\mu+\beta+1;x)
$$
has nonnegative coefficients at all powers of $x$ for all $\mu,\beta\geq{0}$.
Then $\mu\to{f(\mu;x)}$ is  coefficient-wise Wright log-concave on $[0,\infty)$ for $\alpha\in\N$, $\beta\geq\alpha-1$.
\end{lemma}

\begin{theorem}\label{th:psilinearization}
Suppose $\alpha\in\N$. Then the following identity holds true
\begin{multline}\label{eq:psialphabeta}
\frac{\psi_{a,b,\alpha,\beta}(x)}{x^{2a+\alpha+\beta}(1-x)^{2b}}
=\sum\limits_{j=0}^{\alpha-1}\frac{\Gamma(a+b+j)\Gamma(a+b+\beta+\alpha-j-1)}
{[\Gamma(b)]^2\Gamma(a+j+1)\Gamma(a+\beta+\alpha-j)}\times
\\
\left\{\frac{a+b+j}{a+j+1}{_2F_1}\!\!\left(\!\!\begin{array}{l}1,a+b+j+1\\a+j+2\end{array};x\!\right)
-\frac{a+b+\beta+\alpha-1-j}{a+\beta+\alpha-j}{_2F_1}\!\!\left(\!\!\begin{array}{l}1,a+b+\beta+\alpha-j\\a+\beta+\alpha-j+1\end{array};x\!\right)\right\}.
\end{multline}
\end{theorem}
\textbf{Proof.}  First we investigate the function $\psi_{a,b,1,\mu}(x)$ defined in (\ref{eq:psi-def}). Using (\ref{eq:I2F1}) we compute
\begin{multline*}
\frac{\psi_{a,b,1,\mu}(x)[\Gamma(b)]^2}{x^{2a+\mu+1}(1-x)^{2b}}=\overbrace{\frac{\Gamma(a+b)\Gamma(a+b+\mu)}{\Gamma(a+1)\Gamma(a+\mu+1)}}^{=A_{a,b,1,\mu}}
\left\{\frac{a+b}{a+1}{_2F_1}\!\!\left(\!\!\begin{array}{l}1,a+b+1\\a+2\end{array};x\!\right){_2F_1}\!\!\left(\!\!\begin{array}{l}1,a+b+\mu\\a+\mu+1\end{array};x\!\right)\right.
\\
\left.-\frac{a+b+\mu}{a+\mu+1}{_2F_1}\!\!\left(\!\!\begin{array}{l}1,a+b\\a+1\end{array};x\!\right){_2F_1}\!\!\left(\!\!\begin{array}{l}1,a+b+\mu+1\\a+\mu+2\end{array};x\!\right)\right\}
=A_{a,b,1,\mu}\sum\limits_{m=0}^{\infty}x^m
\\
\times\sum\limits_{k=0}^{m}\left\{\frac{(a+b)(a+b+1)_k(a+b+\mu)_{m-k}}{(a+1)(a+2)_{k}(a+\mu+1)_{m-k}}
-\frac{(a+b+\mu)(a+b)_k(a+b+\mu+1)_{m-k}}{(a+\mu+1)(a+1)_{k}(a+\mu+2)_{m-k}}\right\}
\\
=\frac{A_{a,b,1,\mu}(b-1)}{(a+1)(a+\mu+1)}\sum\limits_{m=0}^{\infty}x^m\sum\limits_{k=0}^{m}\frac{(a+b)_k(a+b+\mu)_{m-k}}{(a+2)_{k}(a+\mu+2)_{m-k}}(m-2k+\mu)
\\
=A_{a,b,1,\mu}\left\{\frac{a+b}{a+1}{_2F_1}\!\!\left(\!\!\begin{array}{l}1,a+b+1\\a+2\end{array};x\!\right)
-\frac{a+b+\mu}{a+\mu+1}{_2F_1}\!\!\left(\!\!\begin{array}{l}1,a+b+\mu+1\\a+\mu+2\end{array};x\!\right)\right\},
\end{multline*}
where we have used the summation formula
$$
\sum\limits_{k=0}^{m}\frac{(a+b)_k(a+b+\mu)_{m-k}}{(a+2)_{k}(a+\mu+2)_{m-k}}(m-2k+\mu)
=\frac{(a+1)(a+\mu+1)}{(b-1)}\left[\frac{(a+b)_{m+1}}{(a+1)_{m+1}}-\frac{(a+b+\mu)_{m+1}}{(a+\mu+1)_{m+1}}\right]
$$
demonstrated in \cite[p.72]{KKIZV}.  Next we set $g(a)=I_x(a,b)/x^a(1-x)^b$ and calculate
\begin{multline*}
\frac{\psi_{a,b,\alpha,\beta}}{x^{2a+\alpha+\beta}(1-x)^{2b}}=g(a+\alpha)g(a+\beta)-g(a)g(a+\alpha+\beta)
\\[4pt]
=\!\![g(a+\alpha)g(a+\beta)-g(a+\alpha-1)g(a+\beta+1)]
+[g(a+\alpha-1)g(a+\beta+1)-g(a+\alpha-2)g(a+\beta+2)]
\\[4pt]
+\cdots+[g(a+2)g(a+\beta+\alpha-2)-g(a+1)g(a+\beta+\alpha-1)]+[g(a+1)g(a+\beta+\alpha-1)-g(a)g(a+\beta+\alpha)]
\\
=\frac{1}{[\Gamma(b)]^2}\sum\limits_{j=0}^{\alpha-1}A_{a+j,b,1,\alpha+\beta-2j-1}
\left\{\frac{a+b+j}{a+j+1}{_2F_1}\!\!\left(\!\!\begin{array}{l}1,a+b+j+1\\a+j+2\end{array};x\!\right)\right.
\\[4pt]
\left.-\frac{a+b+\beta+\alpha-1-j}{a+\beta+\alpha-j}{_2F_1}\!\!\left(\!\!\begin{array}{l}1,a+b+\beta+\alpha-j\\a+\beta+\alpha-j+1\end{array};x\!\right)\right\}.
\end{multline*}
Substituting the value of $A_{a+j,b,1,\alpha+\beta-2j-1}$ into this formula yields (\ref{eq:psialphabeta}).~$\hfill\square$

The next corollary refines Theorem~\ref{th:logconc-a} and gives partial solution to Conjecture~\ref{cj:psi}.
\begin{corollary}\label{cr:psi-k}
Suppose $\alpha\in\N$. Then the coefficients $\psi_k$ defined in \emph{(\ref{eq:psi-def})} can be computed by
\begin{multline}\label{eq:psi-k}
\psi_k=
\sum\limits_{j=0}^{\alpha-1}\frac{\Gamma(a+b+j)\Gamma(a+b+\beta+\alpha-j-1)}
{[\Gamma(b)]^2\Gamma(a+j+1)\Gamma(a+\beta+\alpha-j)}
\\
\times\left\{\frac{(a+b+j)_{k+1}}{(a+j+1)_{k+1}}
-\frac{(a+b+\beta+\alpha-1-j)_{k+1}}{(a+\beta+\alpha-j)_{k+1}}\right\}.
\end{multline}
Furthermore,  if $0<b<1$ then $\psi_k<0$ for all $\alpha,\beta>0$  and  if $b>1$ then
$\psi_k>0$ for $\beta\ge\alpha-1$,  $\beta>0$ and $\alpha\in\N$.
\end{corollary}
\textbf{Proof.}  Formula (\ref{eq:psi-k}) follows immediately from (\ref{eq:psialphabeta}).
The case $0<b<1$ falls under conditions of \cite[Theorem~3]{KKJMAA} and the claim follows from this theorem.
Further, according to Lemma~\ref{lm:q-log} for $\beta\ge\alpha-1$ the sign of $\psi_k=\psi_k(\alpha,\beta)$ is determined by the sign of
$\psi_k(1,\beta)$. In view of $(a)_k=\Gamma(a+k)/\Gamma(a)$ for  $\alpha=1$ we have
$$
\psi_{k-1}(1,\beta)=\frac{\Gamma(a+b)\Gamma(a+b+\mu)}{\Gamma(a+1)\Gamma(a+\mu+1)}\left[\frac{\Gamma(a+k+b)\Gamma(a+1)}{\Gamma(a+b)\Gamma(a+k+1)}
-\frac{\Gamma(a+\mu+k+b)\Gamma(a+\mu+1)}{\Gamma(a+\mu+b)\Gamma(a+\mu+k+1)}\right].
$$
Define the function
$$
f(x)=\frac{\Gamma(a+\mu+k+x)\Gamma(a+x)}{\Gamma(a+k+x)\Gamma(a+\mu+x)}.
$$
This function is decreasing according to \cite[Theorem~6]{BusIsm}.  It is straightforward
that $\psi_{k-1}(1,\beta)<0$ is equivalent to $f(1)<f(b)$ while $\psi_{k-1}(1,\beta)>0$ is equivalent to $f(1)>f(b)$,
so that the claim follows from the decrease of $f(x)$.$\hfill\square$

This corollary can also be derived from our previous results in \cite[Theorem~3]{KKJMAA} and \cite[Theorem~2]{KKIZV},
however, the above derivation has the advantage of presenting the explicit formulas for the coefficients $\psi_k$ which
do not follow from these references.

\textbf{Remark.} Substituting (\ref{eq:I2F1}) for $I_x$ into the definition of $\psi_{a,b,\alpha,\beta}$ given in (\ref{eq:psi-def}) and changing notation we can rewrite (\ref{eq:psialphabeta}) as the following linearization identity for the quadratic form in Gauss hypergeometric functions
\begin{multline}\label{eq:2F1-linearization1}
\frac{(a)_{n+1}}{(c)_{n+1}}{_2F_1}\!\!\left(\!\!\begin{array}{l}a+n+1,1\\c+n+1\end{array};x\!\right)
{_2F_1}\!\!\left(\!\!\begin{array}{l}a+\beta,1\\c+\beta\end{array};x\!\right)
\\
-
\frac{(a+\beta)_{n+1}}{(c+\beta)_{n+1}}{_2F_1}\!\!\left(\!\!\begin{array}{l}a,1\\c\end{array};x\!\right)
{_2F_1}\!\!\left(\!\!\begin{array}{l}a+\beta+n+1,1\\c+\beta+n+1\end{array};x\!\right)
\\
=
\sum\limits_{j=0}^{n}\frac{(a)_{j}(a+\beta)_{n-j}}{(c)_j(c+\beta)_{n-j}}
\left\{\frac{a+j}{c+j}{_2F_1}\!\!\left(\!\!\begin{array}{l}a+j+1,1\\c+j+1\end{array};x\!\right)\right.
\\
-
\left.\frac{a+\beta+n-j}{c+\beta+n-j}{_2F_1}\!\!\left(\!\!\begin{array}{l}a+\beta+n+1-j,1\\c+\beta+n+1-j\end{array};x\!\right)
\right\}.
\end{multline}
Here $n=\alpha-1$ is any nonnegative integer, while other parameters as well as $x$ can be arbitrary complex by analytic continuation.  In particular when $n=0$ the above identity reduces to (after some rearrangement)
\begin{equation*}
\frac{a}{c}{_2F_1}\!\!\left(\!\!\begin{array}{l}1,a+1\\c+1\end{array};x\!\right)\left[{_2F_1}\!\!\left(\!\!\begin{array}{l}1,a+\beta\\c+\beta\end{array};x\!\right)-1\right]
=
\frac{a+\beta}{c+\beta}{_2F_1}\!\!\left(\!\!\begin{array}{l}1,a+\beta+1\\c+\beta+1\end{array};x\!\right)\left[{_2F_1}\!\!\left(\!\!\begin{array}{l}1,a\\c\end{array};x\!\right)-1\right]
\end{equation*}

Formula (\ref{eq:2F1-linearization1}) leads to a curious combinatorial identity which is presumably new.
\begin{corollary}\label{cr:combin1}
For arbitrary nonnegative integers $m$, $n$ the following identity holds
\begin{multline*}
\sum\limits_{k=0}^{m}\frac{(a)_k(a+\beta)_{m-k}}{(b)_k(b+\beta)_{m-k}}
\left\{\frac{(a+k)_{n+1}}{(b+k)_{n+1}}-\frac{(a+\beta+m-k)_{n+1}}{(b+\beta+m-k)_{n+1}}\right\}
\\
=\sum\limits_{j=0}^{n}\frac{(a)_j(a+\beta)_{n-j}}{(b)_j(b+\beta)_{n-j}}
\left\{\frac{(a+j)_{m+1}}{(b+j)_{m+1}}-\frac{(a+\beta+n-j)_{m+1}}{(b+\beta+n-j)_{m+1}}\right\}.
\end{multline*}
\end{corollary}
\textbf{Proof}.  The claimed identity follows on equating coefficients at
$x^m$ on both sides of the (\ref{eq:2F1-linearization1}) and using the formula
$(\mu)_{n+1}(\mu+n+1)_k=(\mu)_k(\mu+k)_{n+1}$.~~$\hfill\square$

The next corollary presents two-sided bounds for the Tur\'{a}n determinant formed
by shifts of $I_x(a,b)$ in the first parameter.

\begin{corollary}\label{cr:Turan-a}
Suppose $b>1$, $\nu\in\N$ and $a>\nu-1$. Then
\begin{equation}\label{eq:Turan-a}
mx^a(1-x)^b\leq I_x(a,b)^2-I_x(a-\nu,b)I_x(a+\nu,b)\leq MI_x(a,b)^2
\end{equation}
for all $x\in(0,1)$. Here
\begin{equation*}
\begin{split}
&m=\frac{[\Gamma(a+b)]^2}{[\Gamma(b)]^2[\Gamma(a+1)]^2}-\frac{\Gamma(a+b-\nu)\Gamma(a+b+\nu)}{[\Gamma(b)]^2\Gamma(a-\nu+1)\Gamma(a+\nu+1)},
\\
&M=\frac{(a+1)_{\nu}(a+b-\nu)_{\nu}-(a+1-\nu)_{\nu}(a+b)_{\nu}}{(a+1)_{\nu}(a+b-\nu)_{\nu}}.
\end{split}
\end{equation*}
\end{corollary}
\textbf{Proof.} Inequality from below is just a rewriting of the inequality $\psi_{a-\nu,b,\nu,\nu}\ge{x^{2a}(1-x)^{2b}}\psi_0$
which follows from nonnegativity of the coefficients $\psi_k$ for all $k=0,1,\ldots$  To prove the upper bound consider
the function
$$
a\to G(a)=\frac{I_x(a,b)\Gamma(a+1)}{x^{a}(1-x)^{b}\Gamma(a+b)}
=\frac{1}{\Gamma(b)}\sum\limits_{n=0}^{\infty}\frac{(a+b)_n}{(a+1)_{n}}x^n.
$$
According to \cite[Theorem~2]{KSJMAA2010} the above function is log-convex for $b>1$,
 so that for $a>\nu-1$ we have $G(a-\nu)G(a+\nu)\ge[G(a)]^2$.
This  inequality can be rewritten as
$$
I_x(a,b-\nu)I_x(a,b+\nu)\ge\frac{[\Gamma(a+1)]^2\Gamma(a+b+\nu)\Gamma(a+b-\nu)}{[\Gamma(a+b)]^2\Gamma(a+\nu+1)\Gamma(a-\nu+1)}[I_x(a,b)]^2.
$$
Simple rearrangement of the above inequality gives the upper bound in (\ref{eq:Turan-a}). $\hfill\square$

To investigate the power series coefficients of the function $\phi_{a,b,\alpha,\beta}(x)$ defined in (\ref{eq:Wright})
we need the following lemma.
\begin{lemma}\label{lm:Gosper-b}
The following identity holds:
$$
\sum\limits_{k=0}^{m}\frac{(a+b)_k(a+b+\beta)_{m-k}}{(a+1)_k(a+1)_{m-k}}(b(2k-m)+\beta(a+k))
=\frac{a((a+b+\beta)_{m+1}-(a+b)_{m+1})}{(a+1)_m}.
$$
\end{lemma}
\textbf{Proof.} Denote
$$
u_k=\frac{(a+b)_k(a+b+\beta)_{m-k}}{(a+1)_k(a+1)_{m-k}}(b(2k-m)+\beta(a+k)),~~~k=0,1,\ldots,m.
$$
An application of Gosper's algorithm yields anti-differences:
$$
u_k=v_{k+1}-v_{k},~~~v_k=-\frac{a^2(a+b)_k(a+b+\beta)_{m-k+1}}{(a)_{k}(a)_{m-k+1}},~~~k=0,1,\ldots,m+1.
$$
Straightforward substitution provides verification of this formula. Hence, we get
$$
\sum\limits_{k=0}^{m}u_k=v_{m+1}-v_0=
\frac{a((a+b+\beta)_{m+1}-(a+b)_{m+1})}{(a+1)_m}.~~~~~~\square
$$

\begin{theorem}\label{th:philinearization}
Suppose $\alpha\in\N$. Then the following identity holds true
\begin{multline}\label{eq:phialphabeta}
\frac{\phi_{a,b,\alpha,\beta}(x)}{x^{2a}(1-x)^{2b+\alpha+\beta}}
=\sum\limits_{j=0}^{\alpha-1}\frac{a\Gamma(a+b+j)\Gamma(a+b+\beta+\alpha-j-1)}
{[\Gamma(a+1)]^2\Gamma(b+j+1)\Gamma(b+\beta+\alpha-j)}\times
\\
\left\{(a+b+\beta+\alpha-j-1){_2F_1}\!\!\left(\!\!\begin{array}{l}1,a+b+\beta+\alpha-j\\a+1\end{array};x\!\right)
-(a+b+j){_2F_1}\!\!\left(\!\!\begin{array}{l}1,a+b+j+1\\a+1\end{array};x\!\right)\right\}.
\end{multline}
\end{theorem}
\textbf{Proof.}  First we investigate the function $\phi_{a,b,1,\beta}(x)$ defined in (\ref{eq:Wright}).  We have
\begin{multline*}
\frac{\phi_{a,b,1,\beta}(x)}{x^{2a}(1-x)^{2b+1+\beta}}=\sum\limits_{m=0}^{\infty}x^m\sum\limits_{k=0}^{m}
\biggl\{\frac{\Gamma(a+b+1+k)\Gamma(a+b+\beta+m-k)}{\Gamma(b+1)\Gamma(a+1+k)\Gamma(b+\beta)\Gamma(a+1+m-k)}
\\
-\frac{\Gamma(a+b+k)\Gamma(a+b+\beta+1+m-k)}{\Gamma(b)\Gamma(a+1+k)\Gamma(b+1+\beta)\Gamma(a+1+m-k)}\biggr\}
=\sum\limits_{m=0}^{\infty}x^m
\\
\times\sum\limits_{k=0}^{m}\frac{\Gamma(a+b+k)\Gamma(a+b+\beta+m-k)}{\Gamma(b)\Gamma(b+\beta)\Gamma(a+1+k)\Gamma(a+1+m-k)}
\biggl\{\frac{a+b+k}{b}-\frac{a+b+\beta+m-k}{b+\beta}\biggr\}
\\
=\frac{\Gamma(a+b)\Gamma(a+b+\beta)}{\Gamma(b+1)\Gamma(b+\beta+1)\Gamma(a+1)^2}
\sum\limits_{m=0}^{\infty}x^m\sum\limits_{k=0}^{m}\frac{(a+b)_k(a+b+\beta)_{m-k}}{(a+1)_k(a+1)_{m-k}}(b(2k-m)+\beta(a+k)).
\end{multline*}
Application of Lemma~\ref{lm:Gosper-b} then gives
\begin{multline}\label{eq:phi1mu}
\frac{\phi_{a,b,1,\beta}}{x^{2a}(1-x)^{2b+1+\beta}}=\frac{a\Gamma(a+b)\Gamma(a+b+\beta)}{\Gamma(b+1)\Gamma(b+\beta+1)\Gamma(a+1)^2}
\sum\limits_{m=0}^{\infty}\frac{(a+b+\beta)_{m+1}-(a+b)_{m+1}}{(a+1)_m}x^m
\\
=\frac{a\Gamma(a+b)\Gamma(a+b+\beta)}{\Gamma(b+1)\Gamma(b+\beta+1)\Gamma(a+1)^2}
\left\{(a+b+\beta){_2F_1}\!\!\left(\!\!\begin{array}{l}1,a+b+\beta+1\\a+1\end{array};x\!\right)\right.
\\
\left.-(a+b){_2F_1}\!\!\left(\!\!\begin{array}{l}1,a+b+1\\a+1\end{array};x\!\right)\right\}.
\end{multline}
Next denote $f(b)=I_x(a,b)/x^a(1-x)^b$.  For $\alpha\in\N$ we have
\begin{multline*}
\frac{\phi_{a,b,\alpha,\beta}}{x^{2a}(1-x)^{2b+\alpha+\beta}}=f(b+\alpha)f(b+\beta)-f(b)f(b+\alpha+\beta)=[f(b+\alpha)f(b+\beta)
\\
-f(b+\alpha-1)f(b+\beta+1)]+[f(b+\alpha-1)f(b+\beta+1)-f(b+\alpha-2)f(b+\beta+2)]+\cdots
\\
+[f(b+2)f(b+\beta+\alpha-2)-f(b+1)f(b+\beta+\alpha-1)]+[f(b+1)f(b+\beta+\alpha-1)-f(b)f(b+\beta+\alpha)].
\end{multline*}
Each expression in brackets has the form  $x^{-2a}(1-x)^{-2b-\alpha-\beta}\phi_{a,b+\alpha-j,1,\beta-\alpha+2j-1}$, $j=1,2,\ldots,\alpha$.
Hence, we can apply formula (\ref{eq:phi1mu}) to each such bracket which yields (\ref{eq:phialphabeta}).~$\hfill\square$

The next corollary refines Theorem~\ref{th:logconc-b} and gives partial solution to Conjecture~\ref{cj:phi}.
\begin{corollary}\label{cr:psi-k}
Suppose $\alpha\in\N$. Then the coefficients $\phi_k$ defined in \emph{(\ref{eq:phi-expanded})} can be computed by
\begin{multline}\label{eq:phi-k}
\phi_k=
\sum\limits_{j=0}^{\alpha-1}\frac{a\Gamma(a+b+j)\Gamma(a+b+\beta+\alpha-j-1)}
{[\Gamma(a+1)]^2\Gamma(b+j+1)\Gamma(b+\beta+\alpha-j)(a+1)_k}\times
\\
\left\{(a+b+\beta+\alpha-j-1)_{k+1}-(a+b+j)_{k+1}\right\}.
\end{multline}
If, furthermore, $\beta\ge\alpha-1$ and $\beta>0$ then $\phi_k>0$.
\end{corollary}
\textbf{Proof.} Formula (\ref{eq:phi-k}) follows immediately from (\ref{eq:phialphabeta}).
According to Lemma~\ref{lm:q-log} for $\beta\ge\alpha-1$ the sign of $\phi_k=\phi_k(\alpha,\beta)$ is determined by the sign of
$\phi_k(1,\beta)$. In view of (\ref{eq:phi1mu}) for  $\alpha=1$ we have
$$
\phi_k(1,\beta)=\frac{a\Gamma(a+b)\Gamma(a+b+\beta)}{\Gamma(b+1)\Gamma(b+\beta+1)\Gamma(a+1)^2(a+1)_k}
((a+b+\beta)_{k+1}-(a+b)_{k+1}),
$$
which is clearly positive for $\beta>0$. $\hfill\square$

\textbf{Remark}. Substituting (\ref{eq:I2F1}) for $I_x$ into the definition of $\phi_{a,b,\alpha,\beta}$ given in (\ref{eq:Wright}) and changing notation we can rewrite (\ref{eq:phialphabeta}) as the following linearization identity for the quadratic form in Gauss hypergeometric functions
\begin{multline}\label{eq:2F1linearization2}
\frac{(a+b)_{n+1}}{(b)_{n+1}}{_2F_1}\!\!\left(\!\!\begin{array}{l}a+b+n+1,1\\a+1\end{array};x\!\right)
{_2F_1}\!\!\left(\!\!\begin{array}{l}a+b+\beta,1\\a+1\end{array};x\!\right)
\\
-
\frac{(a+b+\beta)_{n+1}}{(b+\beta)_{n+1}}{_2F_1}\!\!\left(\!\!\begin{array}{l}a+b,1\\a+1\end{array};x\!\right)
{_2F_1}\!\!\left(\!\!\begin{array}{l}a+b+\beta+n+1,1\\a+1\end{array};x\!\right)
\\
=\frac{a}{b(b+\beta)}\sum\limits_{j=0}^{n}\frac{(a+b)_{j}(a+b+\beta)_{n-j}}{(b+1)_{j}(b+\beta+1)_{n-j}}\times
\\
\left\{\!(a+b+\beta+n-j){_2F_1}\!\!\left(\!\!\!\begin{array}{l}a+b+\beta+n+1-j,1\\a+1\end{array}\!\!;x\!\right)\!
-\!(a+b+j){_2F_1}\!\!\left(\!\!\!\begin{array}{l}a+b+j+1,1\\a+1\end{array}\!\!;x\!\right)\!\right\}.
\end{multline}
Here $n=\alpha-1$ is any nonnegative integer, while other parameters as well as $x$ can
be arbitrary complex by analytic continuation.  In particular, when $n=0$ the above identity
after some rearrangement reduces to
\begin{multline*}
\frac{a}{a-c+1}{_2F_1}\!\!\left(\!\!\begin{array}{l}1,a+1\\c\end{array};x\!\right)
\left[{_2F_1}\!\!\left(\!\!\begin{array}{l}1,b\\c\end{array};x\!\right)
+\frac{c-1}{b-c+1}\right]
\\
=
\frac{b}{b-c+1}{_2F_1}\!\!\left(\!\!\begin{array}{l}1,b+1\\c\end{array};x\!\right)
\left[{_2F_1}\!\!\left(\!\!\begin{array}{l}1,a\\c\end{array};x\!\right)
+\frac{c-1}{a-c+1}\right].
\end{multline*}

Formula (\ref{eq:2F1linearization2}) leads to a curious combinatorial identity which is presumably new.
\begin{corollary}\label{cr:combinat2}
For arbitrary nonnegative integers $m$, $n$ the following identity holds
\begin{multline*}
\frac{a}{b(b+\beta)}\sum\limits_{j=0}^{n}\frac{(a+b)_j(a+b+\beta)_{n-j}}{(b+1)_{j}(b+1+\beta)_{n-j}}
\{(a+b+\beta+n-j)_{m+1}-(a+b+j)_{m+1}\}
\\
=(a+1)_m\sum\limits_{k=0}^{m}\frac{(a+b)_k(a+b+\beta)_{m-k}}{(a+1)_{k}(a+1)_{m-k}}
\left\{\frac{(a+b+k)_{n+1}}{(b)_{n+1}}-\frac{(a+b+\beta+m-k)_{n+1}}{(b+\beta)_{n+1}}\right\}.
\end{multline*}
\end{corollary}
\textbf{Proof}.  The claimed identity follows on equating coefficients at
$x^m$ on both sides of the (\ref{eq:2F1linearization2}) and using the formula
$(\mu+n+1)_k(\mu)_{n+1}=(\mu)_k(\mu+k)_{n+1}$. $\hfill\square$

The next corollary presents two-sided bounds for the Tur\'{a}n determinant formed
by shifts of $I_x(a,b)$ in the second parameter.

\begin{corollary}\label{cr:Turan-b}
Suppose $\nu$ is a positive integer $b>\nu$, $a>0$. Then
\begin{equation}\label{eq:Turan-b}
mx^{2a}(1-x)^{2b}\le [I_x(a,b)]^2-I_x(a,b+\nu)I_x(a,b-\nu)\le M[I_x(a,b)]^2
\end{equation}
for all $x\in(0,1)$.  Here
\begin{equation*}\label{eq:mM}
\begin{split}
&m=\frac{[\Gamma(a+b)]^2}{[\Gamma(a+1)]^2[\Gamma(b)]^2}-\frac{\Gamma(a+b-\nu)\Gamma(a+b+\nu)}{[\Gamma(a+1)]^2\Gamma(b-\nu)\Gamma(b+\nu)},
\\[5pt]
&M=\frac{(b)_{\nu}(a+b-\nu)_{\nu}-(a+b)_{\nu}(b-\nu)_{\nu}}{(b)_{\nu}(a+b-\nu)_{\nu}}.
\end{split}
\end{equation*}
\end{corollary}
\textbf{Proof.} Inequality from below is just a rewriting of the inequality $\phi_{a,b-\nu,\nu,\nu}\ge{x^{2a}(1-x)^{2b}}\phi_0$
which follows from nonnegativity of the coefficients $\phi_k$ for all $k=0,1,\ldots$  To prove the upper bound consider
the function
$$
b\to F(b)=\frac{I_x(a,b)\Gamma(b)}{x^{a}(1-x)^{b}\Gamma(a+b)}=\sum\limits_{n=0}^{\infty}\frac{(a+b)_n}{\Gamma(a+n+1)}x^n
=\sum\limits_{n=0}^{\infty}\frac{(a+b)_nn!}{\Gamma(a+n+1)}\frac{x^n}{n!}.
$$
It is easy to check that the coefficients $\{n!/\Gamma(a+n+1)\}_{n=0}^{\infty}$ form a log-convex sequence. Then according to
\cite[Theorem~1]{KSJMAA2010} the above function is log-convex, so that for $b>\nu$ we have $F(b-\nu)F(b+\nu)\ge[F(b)]^2$.
This  inequality can be rewritten as
$$
I_x(a,b-\nu)I_x(a,b+\nu)\ge\frac{[\Gamma(b)]^2\Gamma(a+b+\nu)\Gamma(a+b-\nu)}{[\Gamma(a+b)]^2\Gamma(b+\nu)\Gamma(b-\nu)}[I_x(a,b)]^2.
$$
Simple rearrangement of the above inequality gives the upper bound in (\ref{eq:Turan-b}). $\hfill\square$

\paragraph{4. Acknowledgements.} I am indebted to Professor Skip Garibaldi
(Institute for Pure and Applied Mathematics at UCLA) for sharing the log-concavity
problem with me and for explaining the motivation.  This research was supported by
the Ministry of Education of the Russian Federation (project 1398.2014) and by Russian Foundation for Basic
Research (project 15-56-53032).

\end{document}